\documentclass{amsart}
\usepackage{latexsym}
\usepackage{amsmath,amsthm,amsfonts,amscd,eucal}
\usepackage{graphicx}
\usepackage{epsfig}

\numberwithin{equation}{section}

\def\ca{{\mathcal A}}

\def\cc{{\mathcal C}}

\def\cf{{\mathcal F}}

\def\cn{{\mathcal N}}

\def\cp{{\mathcal P}}

\def\car{{\mathcal R}}

\def\bc{{\mathbb C}}

\def\bn{{\mathbb N}}
\def\br{{\mathbb R}}
\def\bz{{\mathbb Z}}

\def\a{\alpha}

\def\g{\gamma}        \def\G{\Gamma}
        \def\D{\Delta}
\def\eps{\varepsilon}

\def\th{\vartheta}

\def\l{\lambda}       \def\La{\Lambda}

\def\n{\nu}

\def\r{\rho}
\def\s{\sigma}        
     
\def\t{\tau}

        \def\O{\Omega}

\newcommand{\diag}{\text{diag}}

\newcommand{\Det}{\text{det}}
\newcommand{\set}[1]{\left\{#1\right\}}
\def\itm#1{\item{$(#1)$}}

\newcommand{\conv}{\text{conv}\,}
\def\ov{\overline}

\newcommand{\Ci}{\cc}  
\newcommand{\Nt}{\cc^{\text{notail}}}  
\newcommand{\Ta}{\cc^{\text{tail}}}  
\renewcommand{\Re}{\car}  
\renewcommand{\Pr}{\cp}  

\newtheorem{Thm}{Theorem}[section]
\newtheorem{Cor}[Thm]{Corollary}
\newtheorem{Prop}[Thm]{Proposition}
\newtheorem{Lemma}[Thm]{Lemma}

\theoremstyle{definition}
\newtheorem{Dfn}[Thm]{Definition}
\newtheorem{exmp}[Thm]{Example}

\theoremstyle{remark}
\newtheorem{rem}[Thm]{Remark} 
\newtheorem{ack}{Acknowledgement} 

\begin{document}

\title[Ihara zeta function for periodic simple graphs]{Ihara zeta
functions for periodic\\ simple graphs}
\author{Daniele Guido, Tommaso Isola, Michel L. Lapidus}%
\date{August 9, 2006}%
\address{(D.G., T.I.) Dipartimento di Matematica, Universit\`a di
Roma ``Tor
Vergata'', I--00133 Roma, Italy.}%
\email{guido@mat.uniroma2.it, isola@mat.uniroma2.it}%
\address{(M.L.L.) Department of Mathematics, University of California,
Riverside, CA 92521-0135, USA.}%
\email{lapidus@math.ucr.edu}%
\thanks{The first and second authors were partially
supported by MIUR, GNAMPA and by
the European Network ``Quantum Spaces - Noncommutative Geometry"
HPRN-CT-2002-00280. The third author was partially supported
by the National Science Foundation, the Academic Senate of the 
University of California, and GNAMPA}%
\subjclass[2000]{05C25; 05C38; 46Lxx; 11M41. }%
\keywords{Periodic graphs, Ihara zeta function,
 analytic determinant, determinant formula, functional equations.}%

\maketitle
\bigskip

\begin{abstract}
 The definition and main properties of the Ihara zeta function for
 graphs are reviewed, focusing mainly on the case of periodic simple
 graphs.  Moreover, we give a new proof of the associated determinant
 formula, based on the treatment developed by Stark and Terras for
 finite graphs. 
\end{abstract}



\section{Introduction}

 The zeta functions associated to finite graphs by Ihara \cite{Ihara},
 Hashimoto \cite{HaHo,Hashi}, Bass \cite{Bass} and others, combine
 features of Riemann's zeta function, Artin L-functions, and Selberg's
 zeta function, and may be viewed as analogues of the Dedekind zeta
 functions of a number field.  They are defined by an Euler product
 and have an analytic continuation to a meromorphic function
 satisfying a functional equation.  They can be expressed as the
 determinant of a perturbation of the graph Laplacian and, for
 Ramanujan graphs, satisfy the Riemann hypothesis \cite{StTe}.
 
 The first attempt in this context to study infinite graphs was made
 by Grigorchuk and $\dot{\text{Z}}$uk \cite{GrZu}, who considered
 graphs obtained as a suitable limit of a sequence of finite graphs.
 They proved that their definition does not depend on the
 approximating sequence in case of Cayley graphs of finitely generated
 residually finite groups, and, more generally, in case of graphs
 obtained as Schreier graphs of a pair $(G,H)$ of a finitely generated
 group $G$ and a separable subgroup $H$.
 
 The definition of the zeta function was extended to (countable)
 periodic graphs  by Clair and Mokhtari-Sharghi in
 \cite{ClMS1}, where the determinant formula  has been proved. 
 They deduce this result as a specialization of the treatment of 
 group actions on trees (the so-called theory of tree lattices, as 
 developed by Bass, Lubotzky and others, see \cite{BaLu}).
 
 The purpose of this work is to give a more direct proof of that
 result, for the case of periodic simple graphs with a free action. 
 We hope that our treatment, being quite elementary, could be useful
 for someone seeking an introduction to the subject.  In a sequel to
 this paper \cite{GILa03}, we shall prove that for periodic amenable
 graphs, the Ihara zeta function can be approximated by the zeta
 functions of a suitable sequence of finite graphs, thereby answering
 in the affirmative a question raised by Grigorchuk and
 $\dot{\text{Z}}$uk in \cite{GrZu}.
  
 In order to provide a self-contained approach to the subject, we
 start by recalling the definition and some properties of the zeta
 function for finite graphs.  Then, after having introduced some
 preliminary notions, we define in Section \ref{sec:Zeta} the analogue
 of the Ihara zeta function, and show that it is a holomorphic
 function, while, in Section \ref{sec:DetFormula}, we prove a
 corresponding determinant formula.  The latter requires some care,
 because it involves the definition and properties of a determinant
 for bounded operators (acting on an infinite dimensional Hilbert
 space and) belonging to a von Neumann algebra with a finite trace. 
 This question is addressed in Section \ref{sec:AnalyticDet}.  In the
 final section, we establish several functional equations.
  
 In closing this introduction, we note that the operator-algebraic
 techniques used here are introduced by the authors in \cite{GILa01}
 in order to study the Ihara zeta functions attached to a new class of
 infinite graphs, called self-similar fractal graphs.
  

  \section{Zeta function for finite graphs}
 
 The Ihara zeta function is defined by means of equivalence classes of
 prime cycles.  Therefore, we need to introduce some terminology from
 graph theory, following \cite{Serre, StTe} with some modifications. 
 
 A {\it graph} $X=(VX,EX)$ consists of a collection $VX$ of objects,
 called {\it vertices}, and a collection $EX$ of objects called
 (oriented) {\it edges}, together with two maps $e\in EX\mapsto
 (o(e),t(e))\in VX\times VX$ and $e\in EX\mapsto \ov{e}\in EX$,
 satisfying the following conditions: $ \ov{\ov{e}}=e$,
 $o(\ov{e})=t(e)$, $\forall e\in EX$.  The vertex $o(e)$ is called the
 {\it origin} of $e$, while $t(e)$ is called the {\it terminus} of
 $e$.  The couple $\set{e,\ov{e}}$ is called a {\it geometric edge}.
 A graph is called {\it simple} if $EX \subset \set{(u,v)\in VX\times
 VX: u \neq v}$, $o(u,v) = u$, $t(u,v) = v$, $\ov{(u,v)} = (v,u)$;
 therefore, the set of geometric edges can be identified with a set of
 unordered pairs of distinct vertices.  Observe that in the literature
 what we have called graph is also called a multigraph, while a simple
 graph is also called a graph.  We will only deal with simple graphs.
 The edge $e=\set{u,v}$ is said to join the vertices $u,v$, while $u$
 and $v$ are said to be {\it adjacent}, which is denoted $u\sim v$.  A
 {\it path} (of length $m$) in $X$ from $v_0\in VX$ to $v_m\in VX$, is
 $(v_{0},\ldots,v_{m})$, where $v_{i}\in VX$ and $v_{i+1}\sim v_i$,
 for $i=0,...,m-1$.  In the following, we denote by $|C|$ the length
 of a path $C$.  A path is {\it closed} if $v_{m}=v_{0}$.  A graph is
 {\it connected} if there is a path between any pair of distinct
 vertices.
   
 \begin{Dfn}[Proper closed Paths]
    \label{def:redPath} 
    
    \itm{i} A path in $X$ has {\it backtracking} if $v_{i-1}=v_{i+1}$,
    for some $i\in\{1,\ldots,m-1\}$.  A path with no backtracking is
    also called {\it proper}.  Denote by $\Ci$ the set of proper
    closed paths.
       
    \itm{ii} A proper closed path $C=(v_{0},\ldots,v_{m}=v_{0})$ has a
    {\it tail} if there is $k\in\{1,\ldots,[m/2]-1\}$ s.t.
    $v_{j}=v_{m-j}$, for $j=1,\ldots,k$.  Denote by $\Ta$ the set of
    proper closed paths with tail, and by $\Nt$ the set of proper
    tail-less closed paths, also called {\it reduced} closed paths.
    Observe that $\Ci=\Ta\cup\Nt$, $\Ta\cap\Nt=\emptyset$.
    
    \itm{iii} A reduced closed path is {\it primitive} if it is not
    obtained by going $n\geq 2$ times around some other closed path.
 \end{Dfn}

 \begin{exmp}
     Some examples of non reduced closed paths are shown in figures
     \ref{fig:Backtracking}, \ref{fig:Tail}.
     \begin{figure}[ht]
	 \centering
	 \psfig{file=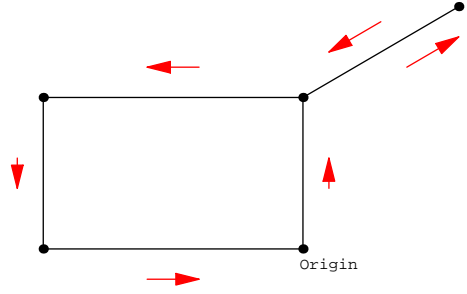,height=1.5in}
	 \caption{Closed path with backtracking}
	 \label{fig:Backtracking}
     \end{figure}
     \begin{figure}[ht]
	 \centering
	 \psfig{file=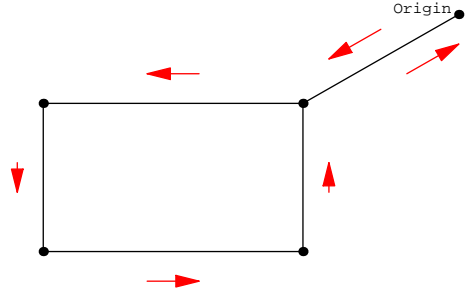,height=1.5in}
	 \caption{Closed path with tail}
	 \label{fig:Tail}
     \end{figure}
 \end{exmp}
 
 We also need an equivalence relation for closed paths
 
  \begin{Dfn}[Cycles]
      Given closed paths $C=(v_{0},\ldots,v_{m}=v_{0})$,
      $D=(w_{0},\ldots,w_{m}=w_{0})$, we say that $C$ and $D$ are {\it
      equivalent}, and write $C\sim_{o} D$, if there is $k$ s.t.
      $w_{j}=v_{j+k}$, for all $j$, where the addition is taken mod
      $m$, that is, the origin of $D$ is shifted $k$ steps w.r.t. the
      origin of $C$.  The equivalence class of $C$ is denoted $[C]_o$.
      An equivalence class is also called a {\it cycle}.  Therefore, a
      closed path is just a cycle with a specified origin.

      Denote by $\Re$ the set of reduced cycles, and by
      $\Pr\subset\Re$ the subset of primitive reduced cycles, also
      called {\it prime} cycles.
 \end{Dfn}
 
 Then Ihara \cite{Ihara} defined the zeta function of a finite graph,
 that is, a graph $X=(VX,EX)$ with $VX$ and $EX$ finite sets, as
 
  \begin{Dfn}[Zeta function]
     $$
     Z_{X}(u) := \prod_{C\in \Pr}
     (1-u^{|C|})^{-1}, \qquad u\in\bc.
     $$
 \end{Dfn}

 Ihara also proved the main result of this theory, though in the
 particular case of regular graphs; subsequently, through the efforts
 of Sunada \cite{Sunada}, Hashimoto \cite{HaHo,Hashi} and Bass
 \cite{Bass}, that result was proved in full generality.  Nowadays,
 there exist many different proofs of Theorem \ref{Thm:Ihara}, $e.g.$
 \cite{StTe,FoZe,KoSu}.  To state it, we need to introduce some more
 notation.  Let us denote by $A=[A(v,w)]$, $v,w\in VX$, the adjacency
 matrix of $X$, that is,
 $$
 A(v,w)=
 \begin{cases} 
     1&\set{v,w}\in EX\\
     0&\text{otherwise.} 
 \end{cases} 
 $$  
 Let $Q:= \diag(\deg(v_{1})-1,\deg(v_{2})-1,\ldots)$, where $\deg(v)$
 is the number of vertices adjacent to $v$, and $\D(u):=I-Au+Qu^{2}$,
 $u\in\bc$, a deformation of the usual Laplacian on the graph, which
 is $\D(1) = (Q+I)-A$.  Then, with $d:= \max_{v\in VX} \deg(v)$, and
 $\chi(X)=|VX|-|EX|$, the Euler characteristic of $X$, we get

 \begin{Thm}[Determinant formula] \label{Thm:Ihara}{\rm
 \cite{Ihara,Sunada,HaHo,Hashi,Bass}}
     $$
     \frac{1}{Z_{X}(u)} = (1-u^{2})^{-\chi(X)}
     \Det(\D(u)), \ \text{ for } |u|<\frac{1}{d-1}.
     $$ 
 \end{Thm}
 
 \begin{exmp}
     We can compute the zeta function of the example shown in figure
     \ref{fig:Example1} by using the determinant formula.  We obtain
     $Z_X(u)^{-1} = (1 - u^2)^2(1 - u)(1 - 2 u)(1 + u + 2 u^2)^3$.
     \begin{figure}[ht]
	 \centering
	 \psfig{file=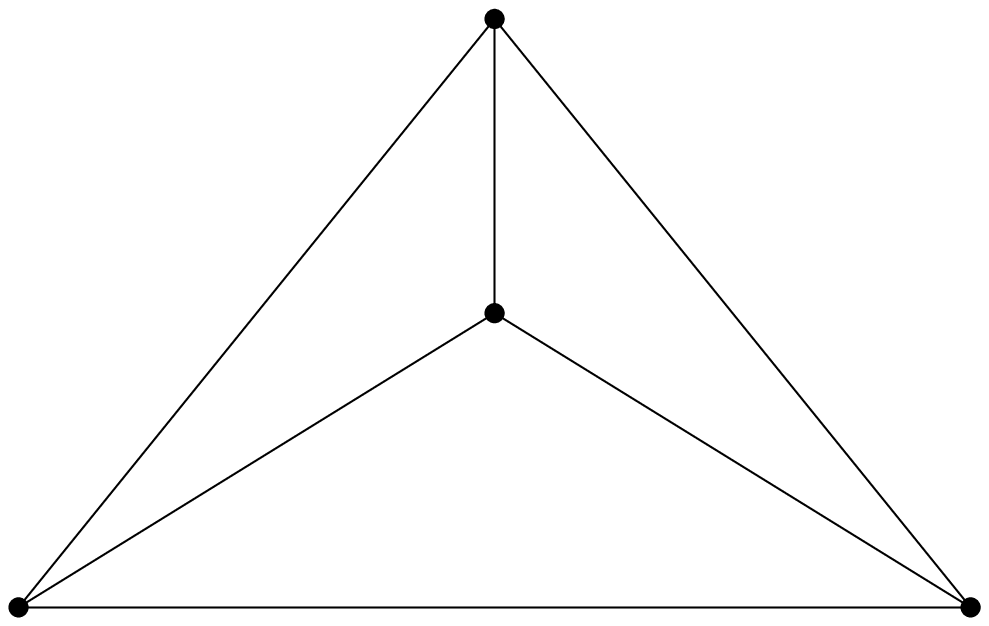,height=1.5in}
	 \caption{A graph}
	 \label{fig:Example1}
     \end{figure}
 \end{exmp}

 The zeta function has been used to establish some properties of the
 graphs.  For example
 
 \begin{Thm} {\rm \cite{Hashi,Hashi1,Bass,North,KoSu}}
     Let $X$ be a finite graph, $r= |EX| - |VX| +1$ the rank of the
     fundamental group $\pi_1(X,x_0)$.  Then $r$ is the order of the
     pole of $Z_X(u)$ at $u=1$.  If $r>1$
     $$ 
     \lim_{u\to 1^-} Z_X(u)(1-u)^r = -\frac{1}{2^r (r-1)\kappa_X},
     $$
     where $\kappa_X$ is the number of spanning trees in $X$.
 \end{Thm}
 
 \begin{Thm} {\rm \cite{Hashi2,HSTe} }
     Let $X$ be a finite graph, $R_X$ be the radius of the greatest
     circle of convergence of $Z_X$.  Denote by $\pi_{n}$ the number
     of prime cycles which have length $n$.  If $g.c.d. \{ |C| :
     C\in\Pr \} =1$, then
     $$
     \pi_{n} \sim \frac{R_X^{-n}}{n},\ n\to\infty.
     $$
 \end{Thm}
 
 \begin{Thm}{\rm  \cite{Ihara,Lubo,StTe} }
     Let $X$ be a finite graph which is $(q+1)$-regular, $i.e.$
     $\deg(v) = q+1$ for all $v\in VX$.  Then the following are
     equivalent
     
     \itm{i}
     $$
     (RH)\ \begin{cases} Z_X(q^{-s})^{-1} = 0 & \\
     \Re s \in(0,1) & \implies \Re s = \frac12.
     \end{cases} 
     $$
     
     \itm{ii} $X$ is a Ramanujan graph, $i.e.$
     $\l\in\s(A),| \lambda | < q+1 \implies |\lambda| \leq 2\sqrt{q}$.
 \end{Thm}
 
 More results on the Ihara zeta function are contained in
 \cite{Serre1,StTe2,StTe3} and in various papers by Mizuno and Sato. 
 In closing this section, we mention a generalization of the Ihara
 zeta function recently introduced by Bartholdi \cite{Bartholdi} and
 studied by Mizuno and Sato (see \cite{MiSa} and references therein).

\section{Periodic simple graphs}

 Let $X=(VX,EX)$ be a simple graph, which we assume to be (countable
 and) with bounded degree, $i.e.$ the degree of the vertices is
 uniformly bounded.  Let $\G$ be a countable discrete subgroup of
 automorphisms of $X$, which acts freely on $X$ ($i.e.$ any $\g\in\G$,
 $\g\neq id$ doesn't have fixed points), and with finite quotient
 $B:=X/\G$.  Denote by $\cf\subset VX$ a set of representatives for
 $VX/\G$, the vertices of the quotient graph $B$.  Let us define a
 unitary representation of $\G$ on $\ell^{2}(VX)$ by $(\l(\g)f)(x):=
 f(\g^{-1}x)$, for $\g\in\G$, $f\in\ell^{2}(VX)$, $x\in V(X)$.  Then
 the von Neumann algebra $\cn(X,\G):= \{ \l(\g) : \g\in\G\}'$ of
 bounded operators on $\ell^{2}(VX)$ commuting with the action of $\G$
 inherits a trace given by $Tr_{\G}(T) = \sum_{x\in\cf} T(x,x)$, for
 $T\in\cn(X,\G)$.

 Let us denote by $A$ the adjacency matrix of $X$.  Then (by
 \cite{Mohar}, \cite{MoWo}) $\|A\|\leq d:=\sup_{v\in VX} \deg(v)
 <\infty$, and it is easy to see that $A\in\cn(X,\G)$.

 For any $m\in\bn$, let us denote by $A_{m}(x,y)$ the number of
 proper paths in $X$, of length $m$, with initial vertex $x$ and
 terminal vertex $y$, for $x,y\in VX$.  Then $A_{1}=A$.  Let $A_{0}:=
 I$ and $Q:= \diag(\deg(v_{1})-1,\deg(v_{2})-1,\ldots)$.  Then
 
 \begin{Lemma}\label{lem:Lemma1}
     \itm{i} $A_{2} = A^{2}-Q-I\in\cn(X,\G)$,
     
     \itm{ii} for $m\geq 3$, $A_{m} = A_{m-1}A-A_{m-2}Q\in\cn(X,\G)$,
     
     \itm{iii} let $\a:= \frac{d+\sqrt{d^{2}+4d}}{2}$; then $\|A_{m}\|
     \leq \a^{m}$, for $m\geq0$.
 \end{Lemma}
 \begin{proof}
     $(i)$ if $x = y$ then $A_{2}(x,x)=0$ because there are no proper
     closed paths of length $2$ starting at $x$, whereas $A^{2}(x,x) =
     \deg(x) = (Q+I)(x,x)$, so that
     $A_{2}(x,x)=A^{2}(x,x)-(Q+I)(x,x)$.  If $x\neq y$, then
     $A^{2}(x,y)$ is the number of paths of length $2$ (necessarily
     proper) from $x$ to $y$, so $A_{2}(x,y)=A^{2}(x,y) =
     A^{2}(x,y)-(Q+I)(x,y)$.  \\
     $(ii)$ for $x,y\in VX$, the sum $\sum_{z\in VX}
     A_{m-1}(x,z)A(z,y)$ counts the proper paths of length $m$ from
     $x$ to $y$, plus additional paths formed of a proper path of
     length $m-2$ from $x$ to $y$ followed by a path of length $2$
     from $y$ to $z$ and back; since the path from $x$ to $y$ and then
     to $z$ is a proper path of length $m-1$ (one of those counted by
     $A_{m-1}(x,z))$, $z$ can only be one of the $deg(y)-1=Q(y,y)$
     vertices adjacent to $y$, the last one being on the proper path
     from $x$ to $y$.  Therefore $\sum_{z\in VX} A_{m-1}(x,z)A(z,y) =
     A_{m}(x,y) + A_{m-2}(x,y)Q(y,y)$, and the statement follows.  \\
     $(iii)$ We have $\|A_{1}\|=\|A\|\leq d$, $\|A_{2}\|\leq 
     d^{2}+d$, and $\|A_{m}\|\leq d(\|A_{m-1}\|+\|A_{m-2}\|$, from 
     which the claim follows by induction.
 \end{proof}
 
 Denote by $\Ci_{m}$ the subset of $\Ci$ consisting of the proper
 closed paths of length $m$, and attach a similar meaning to
 $\Ta_{m}$, $\Nt_{m}$, $\Re_{m}$ and $\Pr_{m}$.
  
 \begin{Lemma}\label{lem:countTail}
     Denote by $t_{m}:= \sum_{x\in \cf} |\{C\in\Ta_{m}: C \text{
     starts at } x \} |$, where $|\cdot|$ denotes the cardinality of a
     set.  Then \itm{i} $t_{1}=t_{2}=0$, and, for $m\geq 3$, $t_{m} =
     Tr_{\G}((Q-I)A_{m-2}) + t_{m-2}$,
     
     \itm{ii} $t_{m} =  
     Tr_{\G}\left( (Q-I)\sum_{j=1}^{[\frac{m-1}{2}]} A_{m-2j}  
\right)$.
 \end{Lemma}
 \begin{proof}
     $(i)$
     Indeed, we have
     \begin{align*}
	 t_{m} &= \sum_{x\in \cf} |\{C\in\Ta_{m}: C \text{ starts at } x 
	 \} |\\
	 &= \sum_{x\in \cf} \sum_{y\sim x} |\{C\in\Ta_{m}: C \text{
	 starts at } x \text{ goes to } y \text{ at first step} \} |\\
	 &= \sum_{y\in \cf} \sum_{x\sim y} |\{C\in\Ta_{m}: C \text{
	 starts at } x \text{ goes to } y \text{ at first step} \} |,
     \end{align*}
     where the last equality follows from the fact that the
     cardinalities above are $\G$-invariant, and we can choose
     $\g\in\G$ for which the second vertex $y$ of $\g C$ is in $\cf$. 
     A path $C$ in the last set goes from $x$ to $y$, then over a
     closed path $D$ of length $m-2$, and then back to $x$.  There are
     two kinds of closed paths $D$ at $y$: those with tails and those
     without.  If $D$ has no tail, then there are $Q(y,y)+1$
     possibilities for $x$ to be adjacent to $y$, but $x$ cannot be on
     $D$ (otherwise, $C$ would have backtracking), which leaves
     $Q(y,y)-1$ possibilities.  If $D$ has a tail, $x$ cannot be on
     $D$ (otherwise, $C$ would have backtracking), which leaves
     $Q(y,y)$ possibilities.  Therefore, we get
     \begin{align*}
	 &\sum_{x\sim y} | \{C\in\Ta_{m}: C \text{
	 starts at } x \text{ goes to } y \text{ at first step} \} |\\
	 &= (Q(y,y)-1)\cdot | \{D\in\Nt_{m-2}: D \text{
	 starts at } y  \} | \\
	 &\qquad + Q(y,y)\cdot | \{D\in\Ta_{m-2}: D \text{
	 starts at } y  \} | \\
	 & = (Q(y,y)-1)\cdot | \{D\in\Ci_{m-2}: D \text{
	 starts at } y  \} | \\
	 & \qquad + | \{D\in\Ta_{m-2}: D \text{
	 starts at } y  \} |\, ,
     \end{align*}
     so that
     \begin{align*}
	 t_{m} &= \sum_{y\in \cf} (Q(y,y)-1)\cdot | \{D\in\Ci_{m-2}: D
	 \text{ starts at } y \} | \\
	 & \qquad + \sum_{y\in \cf} | \{D\in\Ta_{m-2}: D \text{
	 starts at } y  \} | \\
	 & = \sum_{y\in \cf} (Q(y,y)-1)A_{m-2}(y,y) +t_{m-2}\\
	 & = Tr_{\G}((Q-I)A_{m-2})+t_{m-2}.
     \end{align*}
     $(ii)$ Follows from $(i)$.
 \end{proof}
 
 We need to introduce an equivalence relation between reduced cycles.

 \begin{Dfn}[Equivalence relation between reduced cycles]
    Given $C$, $D\in\Re$, we say that $C$ and $D$ are $\G$-{\it
equivalent},
    and write $C \sim_{\G} D$, if there is an isomorphism $\g\in \G$
    s.t. $D=\g(C)$.  We denote by $[\Re]_{\G}$ the set of
    $\G$-equivalence classes of reduced cycles, and analogously for
    the subset $\Pr$.
 \end{Dfn}
 
 For the purposes of the next result, for any closed path 
 $D=(v_{0},\ldots,v_{m}=v_{0})$, we also denote $v_{j}$ by $v_{j}(D)$.
 
 Let us now assume that $C$ is a prime cycle of length $m$.  Then the
 {\it stabilizer}  of $C$ in $\G$ is the subgroup $\G_{C}= \{\g\in\G :
 \g(C)=C\}$ or, equivalently, $\g\in\G_{C}$ if there exists
 $p(\g)\in\bz_{m}$ s.t., for any choice of the origin of $C$,
 $v_{j}(\g C)=v_{j-p}(C)$, for any $j$.  Let us observe that $p(\g)$
 is a group homomorphism from $\G_{C}$ to $\bz_{m}$, which is
 injective because $\G$ acts freely.  As a consequence, $|\G_{C}|$
 divides $m$.
 
 \begin{Dfn}
     Let $C\in\Pr$ and define $\displaystyle \n(C) := 
     \frac{|C|}{|\G_{C}|}$. If $C=D^{k}\in\Re$, where $D\in\Pr$, 
     define $\n(C)=\n(D)$. Observe that $\n(C)$ only depends on 
     $[C]_{\G}\in[\Re]_{\G}$.
 \end{Dfn}

 \begin{Lemma}\label{lem:estim.for.N}
     Let us set $N_{m} := \sum_{[C]_{\G}\in[\Re_{m}]_{\G}}\n(C)$.
Then 
     
     \itm{i} $N_{m} = Tr_{\G}(A_{m}) - t_{m}$,
     
     \itm{ii} $N_{m} \leq d(d-1)^{m-1}|\cf|$.
 \end{Lemma}
 \begin{proof}
     $(i)$ Let us assume that $[C]_{\G}$ is an equivalence class of
     prime cycles in $[\cp_{m}]_{\G}$, and consider the set $U$ of all
     primitive closed paths with the origin in $\cf$ and representing
     $[C]_{\G}$.  If $C$ is such a representative, any other
     representative can be obtained in this way: choose $k\in\bz_{m}$,
     let $\g(k)$ be the (unique) element in $\G$ for which
     $\g(k)v_{k}(C)\in\cf$, and define $C_{k}$ as
     $$
     v_{j}(C_{k})=\g(k)v_{j+k}(C),\ j\in\bz_{m}.
     $$
     If we want to count the elements of $U$, we should know how many
     of the elements $C_{k}$ above coincide with $C$.  For this to
     happen, $\g$ should clearly be in the stabilizer of the cycle
     $[C]_{o}$.  Conversely, for any $\g\in\G_{C}$, there exists
     $p=p(\g)\in\bz_{m}$ such that $\g v_{j}(C)=v_{j-p}(C)$, therefore
     $\g=\g(p)$.  As a consequence, $v_{j}(C_{p(\g)})=\g(p)
     v_{j+p}(C)=v_{j}(C)$, so that $C_{p(\g)}=C$.  We have proved that
     the cardinality of $U$ is equal to $\n(C)$.  The proof for a
     non-prime cycle is analogous.  Therefore,
     \begin{align}
	 N_{m}&= \sum_{[C]_{\G}\in[\Re_{m}]_{\G}} | \{
	 D\in\cc^{notail}_{m}: [D]_{o}\sim_{\G} C, v_{0}(D)\in\cf \} |
\notag \\
	 & = | \{ C\in\cc^{notail}_{m}, v_{0}(C)\in\cf \}
	 | \label{eq:Nm}\\
	 &= | \{ C\in\cc_{m}, v_{0}(C)\in\cf \} | - | \{
	 C\in\cc^{tail}_{m}, v_{0}(C)\in\cf \} | \notag\\
	 &= Tr_{\G}(A_{m}) - t_{m}. \notag
     \end{align}
     $(ii)$ Follows from (\ref{eq:Nm}).

 \end{proof}
 
\section{The Zeta function}\label{sec:Zeta}
 
 In this section, we define the Ihara zeta function for a periodic
 graph, and prove that it is a holomorphic function.  

 \begin{Dfn}[Zeta function] We let
     $$
     Z_{X,\G}(u) := \prod_{[C]_{\G}\in [\Pr]_{\G}} (1-u^{|C|})^{
     -\frac{1}{ |\G_{C}| } }, 
     $$
     for all $u\in\bc$ sufficiently small so that the infinite
product 
     converges.
 \end{Dfn}

 \begin{Lemma}\label{lem:power.series}
     \itm{i} $Z(u):=\prod_{[C] \in [\Pr]_{\G}}
     (1-u^{|C|})^{ -\frac{1}{ |\G_{C}| } }$, defines a holomorphic
     function in $\{u\in\bc: |u|<\frac{1}{d-1}\}$,
     
     \itm{ii} $u\frac{Z'(u)}{Z(u)} = \sum_{m=1}^{\infty}
     N_{m}u^{m}$, for $|u| < \frac{1}{d-1}$,
     
     \itm{iii} $Z(u) = \exp\left( \sum_{m=1}^{\infty}
     \frac{N_{m}}{m}u^{m} \right)$ , for $|u| < \frac{1}{d-1}$.
 \end{Lemma}
 \begin{proof}
     Let us observe that, for  $|u|<\frac{1}{d-1}$, 
     \begin{align*}
	 \sum_{m=1}^{\infty} N_{m}  u^{m} & =
	 \sum_{[C]_{\G}\in [\Re]_{\G}}\n(C)\, u^{|C|} \\
	 & =  \sum_{m=1}^{\infty} \sum_{[C]_{\G}\in
	 [\Pr]_{\G}}\frac{|C|}{|\G_{C}|}\,  u^{|C^{m}|} \\
	 &= \sum_{[C]_{\G}\in
	 [\Pr]_{\G}}\frac{1}{|\G_{C}|}\, \sum_{m=1}^{\infty} |C| u^{|C|m} \\
	 &= \sum_{[C]_{\G}\in [\Pr]_{\G}} \frac{1}{|\G_{C}|}\,  u\frac{d}{du}
	 \sum_{m=1}^{\infty} \frac{u^{|C|m}}{m} \\
	 &= -\sum_{[C]_{\G}\in [\Pr]_{\G}} \frac{1}{|\G_{C}|}\, u\frac{d}{du}
	 \log(1-u^{|C|}) \\
	 & = u\frac{d}{du} \log Z(u),
     \end{align*}
     where, in the last equality we used uniform convergence on
     compact subsets of $\set{u\in\bc: |u|<\frac{1}{d-1}}$.  The rest
     of the proof is clear.
 \end{proof}
 
 \begin{exmp}
     Some examples of cycles with different stabilizers are shown in
     figures \ref{fig:cycle1}, \ref{fig:cycle2}.  They refer to the
     graph in figure \ref{fig:Example2} which is the standard lattice
     graph $X=\bz^{2}$ endowed with the action of the group $\G$ which
     is generated by the rotation by $\frac{\pi}{2}$ around the point
     $P$ and the translations by elements $(m,n)\in\bz^{2}$ acting as
     $(m,n)(v_{1},v_{2}):= (v_{1}+2m,v_{2}+2n)$, for
     $v=(v_{1},v_{2})\in VX=\bz^{2}$.
     \begin{figure}[ht]
	 \centering
	 \psfig{file=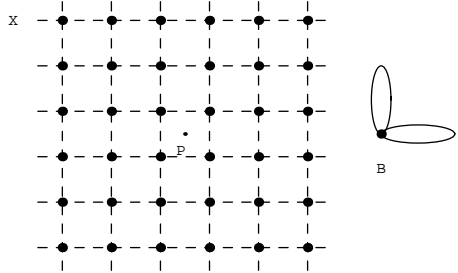,height=1.5in}
	 \caption{A periodic graph $X$ with its quotient $B=X/\G$}
	 \label{fig:Example2}
     \end{figure}
     \begin{figure}[ht]
	 \centering
	 \psfig{file=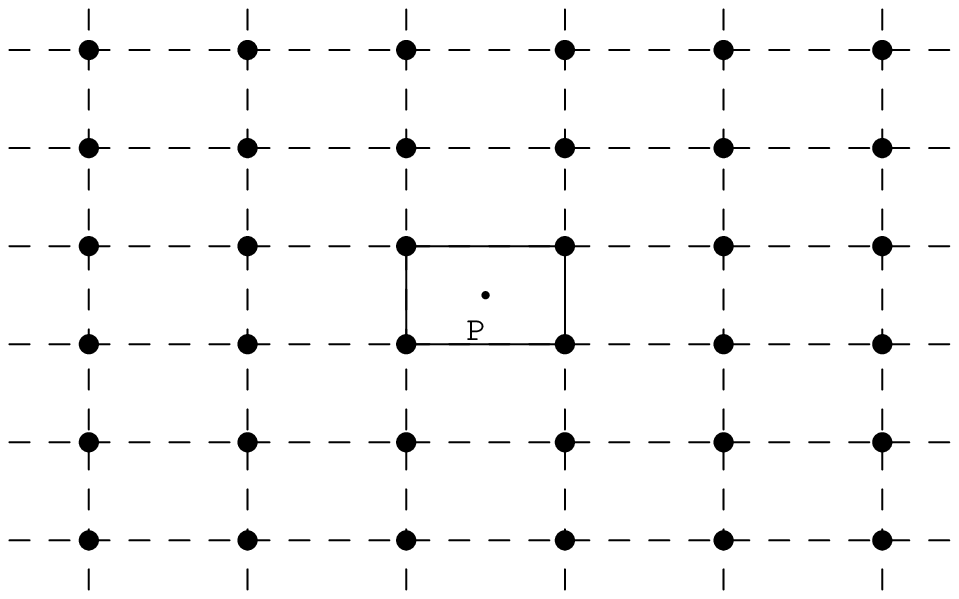,height=1.5in}
	 \caption{A cycle with $|\G_{C}|=4$}
	 \label{fig:cycle1}
     \end{figure}
     \begin{figure}[ht]
	 \centering
	 \psfig{file=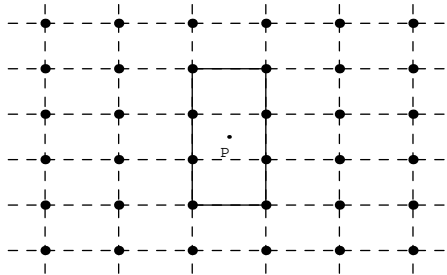,height=1.5in}
	 \caption{A cycle with $|\G_{C}|=2$}
	 \label{fig:cycle2}
     \end{figure}
 \end{exmp}
 
 The interested reader can find the computation of the Ihara zeta
function for several periodic simple graphs in
\cite{GrZu,ClMS1,ClMS2,Clair}.
 
\section[Analytic determinant for von Neumann algebras]{An 
analytic determinant for von Neumann algebras with a finite 
trace}\label{sec:AnalyticDet}

 In this section, we define a determinant for a suitable class of not
 necessarily normal operators in a von Neumann algebra with a finite
 trace.  The results obtained are used in Section \ref{sec:DetFormula}
 to prove a determinant formula for the zeta function.

 In a celebrated paper \cite{FuKa}, Fuglede and Kadison defined a
 positive-valued determinant for finite factors ($i.e.$ von Neumann
 algebras with trivial center and finite trace).  Such a determinant
 is defined on all invertible elements and enjoys the main properties
 of a determinant function, but it is positive-valued.  Indeed, for an
 invertible operator $A$ with polar decomposition $A=UH$, where $U$ is
 a unitary operator and $H:= \sqrt{A^{*}A}$ is a positive self-adjoint
 operator, the Fuglede--Kadison determinant is defined by
 $$
 Det(A)=\exp\, \circ\ \tau\circ\log H,
 $$
 where $\log H$ may be defined via the functional calculus.

 For the purposes of the present paper, we need a determinant which is
 an analytic function.  As we shall see, this can be achieved, but
 corresponds to a restriction of the domain of the determinant
 function and implies the loss of some important properties.  Let
 $(\ca,\tau)$ be a von Neumann algebra endowed with a finite trace. 
 Then, a natural way to obtain an analytic function is to define, for
 $A\in\ca$, $\Det_\t(A)=\exp\, \circ\ \tau\circ\log A$, where
 $$
 \log(A) := \frac{1}{2\pi i} \int_\Gamma \log \lambda (\lambda-A)^{-1}
 d\lambda,
 $$
 and $\Gamma$ is the boundary of a connected, simply connected region
 $\Omega$ containing the spectrum of $A$.  Clearly, once the branch of
 the logarithm is chosen, the integral above does not depend on
 $\Gamma$, provided $\G$ is given as above.

 Then a na\"{\i}ve way of defining $det$ is to allow all elements $A$
 for which there exists an $\Omega$ as above, and a branch of the
 logarithm whose domain contains $\Omega$.  Indeed the following
 holds.

\begin{Lemma}
    Let $A$, $\Omega$, $\Gamma$ be as above, and $\varphi$, $\psi$ two
    branches of the logarithm such that both domains contain $\Omega$.
    Then
    $$
    \exp\, \circ\ \tau\circ\varphi(A) = \exp\, \circ\
    \tau\circ\psi(A).
    $$
\end{Lemma}

\begin{proof}
    The function $\varphi(\lambda)-\psi(\lambda)$ is continuous and
    everywhere defined on $\Gamma$.  Since it takes its values in
    $2\pi i \mathbb{Z}$, it should be constant on $\Gamma$.  Therefore
    \begin{align*}
	\exp\, \circ\ \tau\circ\varphi (A) & = \exp\, \circ\
	\tau\left(\frac{1}{2\pi i} \int_\Gamma 2\pi i n_{0}
	(\lambda-A)^{-1} d\lambda \right) \exp\, \circ\
	\tau\circ\psi(A)\\
	&=\exp\, \circ\ \tau\circ\psi(A).
    \end{align*}
\end{proof}

 The problem with the previous definition is its dependence on the
 choice of $\Omega$.  Indeed, it is easy to see that when
 $A=\begin{pmatrix}1&0\\0&i\end{pmatrix}$ and we choose $\Omega$
 containing $\{e^{i\vartheta},\vartheta\in[0,\pi/2]\}$ and any
 suitable branch of the logarithm, we get $det(A)=e^{i\pi/4}$, by
 using the normalized trace on $2\times 2$ matrices.  On the other
 hand, if we choose $\Omega$ containing
 $\{e^{i\vartheta},\vartheta\in[\pi/2,2\pi]\}$ and a corresponding
 branch of the logarithm, we get $det(A)=e^{5i\pi/4}$.  Therefore, we
 make the following choice.

\begin{Dfn}
    Let $(\ca,\tau)$ be a von Neumann algebra endowed with a finite
    trace, and consider the subset $\ca_{0}=\{A\in\ca : 0\not\in
    \conv\sigma(A)\}$, where $\s(A)$ denotes the spectrum of $A$.  For
    any $A\in\ca_{0}$ we set
    $$
    \Det_\t(A)=\exp\, \circ\ \tau\circ\left(\frac{1}{2\pi i}
    \int_\Gamma \log \lambda (\lambda-A)^{-1} d\lambda\right),
    $$
    where $\Gamma$ is the boundary of a connected, simply connected
    region $\Omega$ containing $\conv\sigma(A)$, and $\log$ is a
    branch of the logarithm whose domain contains $\Omega$.
\end{Dfn}

\begin{Cor}\label{cor:det.analytic}
    The determinant function defined above is well-defined and
    analytic on $\ca_{0}$.
\end{Cor}

We collect some properties of our determinant in the following result.

\begin{Prop}\label{properties}
    Let $(\ca,\tau)$ be a von Neumann algebra endowed with a finite
    trace, $A\in\ca_{0}$.  Then
    
    \item[$(i)$] $\Det_\t(zA)=z^{\t(I)}\Det_\t(A)$, for any
    $z\in\mathbb{C}\setminus\{0\}$,
    
    \item[$(ii)$] if $A$ is normal, and $A=UH$ is its polar
    decomposition, 
    $$\Det_\t (A)=\Det_\t(U)\Det_\t(H),$$
    
    \item[$(iii)$] if $A$ is positive, $\Det_\t(A)=Det(A)$, where the
    latter is the Fuglede-Kadison determinant.
\end{Prop}

\begin{proof}
    $(i)$ If the half-line $\{\rho e^{i\vartheta_0}\in\mathbb{C} :
    \r>0\}$ does not intersect $\conv\sigma(A)$, then the half-line
    $\{\rho e^{i(\vartheta_0+t)}\in\mathbb{C} : \r>0\}$ does not
    intersect $\conv\sigma(zA)$, where $z=re^{it}$.  If $\log$ is the
    branch of the logarithm defined on the complement of the real
    negative half-line, then $\varphi(x)=i(\vartheta_{0}-\pi) +
    \log(e^{-i(\vartheta_{0}-\pi)}x)$ is suitable for defining
    $\Det_\t(A)$, while $\psi(x)=i(\vartheta_{0}+t-\pi) +
    \log(e^{-i(\vartheta_{0}+t-\pi)}x)$ is suitable for defining
    $\Det_\t(zA)$.  Moreover, if $\Gamma$ is the boundary of a
    connected, simply connected region $\Omega$ containing
    $\conv\sigma(A)$, then $z\Gamma$ is the boundary of a connected,
    simply connected region $z\Omega$ containing $\conv\sigma(zA)$. 
    Therefore,
    \begin{align*}
	\Det_\t(zA) &= \exp\, \circ\ \tau\left(\frac{1}{2\pi i}
	\int_{z\Gamma} \psi(\lambda) (\lambda-zA)^{-1}
	d\lambda\right)\\
	&= \exp\, \circ\ \tau\left(\frac{1}{2\pi i} \int_{\Gamma}
	(i(\vartheta_{0}+t-\pi) + \log(e^{-i(\vartheta_{0}+t-\pi)}
	re^{it}\mu)) (\mu-A)^{-1} d\mu\right)\\
	&= \exp\, \circ\ \tau\left((\log r + it)I+\frac{1}{2\pi i}
	\int_{\Gamma} \varphi(\mu) (\mu-A)^{-1} d\mu\right)\\
	&= z^{\t(I)} \Det_\t(A).
    \end{align*}
    $(ii)$ When $A=UH$ is normal, $U=\int_{[0,2\pi]} e^{i\vartheta}\
    du(\th)$, $H=\int_{[0,\infty)}r\ dh(r)$, then $ A =
    \int_{[0,\infty)\times[0,2\pi]} r e^{i\vartheta} \ d(h(r)\otimes
    u(\th))$.  The property $0\not\in\conv\sigma(A)$ is equivalent to
    the fact that the support of the measure $d(h(r)\otimes u(\th))$
    is compactly contained in some open half-plane $$\{\rho
    e^{i\vartheta} : \rho>0, \vartheta \in (\vartheta_{0} - \pi/2,
    \vartheta_{0} +\pi/2)\},$$ or, equivalently, that the support of
    the measure $dh(r)$ is compactly contained in $(0,\infty)$, and
    the support of the measure $d u(\th)$ is compactly contained in
    $(\vartheta_{0} - \pi/2, \vartheta_{0} +\pi/2)$.  Therefore
    $A\in\ca_{0}$ is equivalent to $U,H\in\ca_{0}$.  Then $$\log A =
    \int_{[0,\infty) \times (\vartheta_{0} - \pi/2, \vartheta_{0}
    +\pi/2)} (\log r + i\vartheta) \ d(h(r)\otimes u(\th)),$$ which
    implies that
    \begin{align*}
	\Det_\t(A) &= \exp\, \circ\ \tau\left(\int_{0}^{\infty} \log
	r\ dh(r) + \int_{\vartheta_{0} - \pi/2}^{\vartheta_{0} +\pi/2}
	i\vartheta \ du(\th)\right) \\
	&= \Det_\t(U)\cdot \Det_\t(H).
    \end{align*}
    $(iii)$  Follows by the above argument.
\end{proof}

\begin{rem} 
    We note that the above defined determinant function strongly
    violates the product property $\Det_\t(AB)=\Det_\t(A)\Det_\t(B)$. 
    Indeed, the fact that $A,B\in\ca_{0}$ does not imply
    $AB\in\ca_{0}$, as is seen e.g. by taking
    $A=B=\begin{pmatrix}1&0\\0&i\end{pmatrix}$.  Moreover, even if
    $A,B,AB\in\ca_{0}$ and $A$ and $B$ commute, the product property
    may be violated, as is shown by choosing
    $A=B=\begin{pmatrix}1&0\\0&e^{3i\pi/4}\end{pmatrix}$, and using
    the normalized trace on $2\times 2$ matrices.
\end{rem}
 
\section{The determinant formula}\label{sec:DetFormula}

 In this section, we prove the main result in the theory of Ihara zeta
 functions, which says that $Z$ is the reciprocal of a holomorphic
 function, which, up to a factor, is the determinant of a deformed
 Laplacian on the graph.  We first need some technical results.  Let
 us recall that $d:=\sup_{v\in VX} \deg(v)$, and $\a:=
 \frac{d+\sqrt{d^{2}+4d}}{2}$.
 
 \begin{Lemma}\label{lem:eq.for.A}
     \itm{i} $\left(\sum_{m\geq 0} A_{m}u^{m}\right)(I-Au+Qu^{2}) =
     (1-u^{2})I$, for $|u|<\frac{1}{\a}$, 
     
     \itm{ii} $\left(\sum_{m\geq 0} \left( \sum_{k=0}^{[m/2]}A_{m-2k}
     \right) u^{m}\right)(I-Au+Qu^{2}) = I$, for $|u|<\frac{1}{\a}$.
 \end{Lemma}
 \begin{proof}
     $(i)$ From Lemma \ref{lem:Lemma1} we obtain
     \begin{align*}
	 \biggl(\sum_{m\geq 0} A_{m}u^{m}\biggr)&(I-Au+Qu^{2}) =
	 \sum_{m\geq 0} A_{m}u^{m} - \sum_{m\geq 0}\left(
	 A_{m}Au^{m+1} -A_{m}Qu^{m+2}\right) \\
	 &= \sum_{m\geq 0} A_{m}u^{m} -A_{0}Au -A_{1}Au^{2}
	 +A_{0}Qu^{2} \\
	 & \qquad - \sum_{m\geq 3}\left( A_{m-1}A
	 -A_{m-2}Q\right)u^{m} \\
	 &= \sum_{m\geq 0} A_{m}u^{m} -Au -A^{2}u^{2}
	 +Qu^{2} - \sum_{m\geq 3} A_{m}u^{m} \\
	 &= I +Au +A_{2}u^{2} -Au -A^{2}u^{2}
	 +Qu^{2} \\
	 & = (1-u^{2})I.
     \end{align*}
     $(ii)$ 
     \begin{align*}
	 I &= (1-u^{2})^{-1} \biggl(\sum_{m\geq 0} 
	 A_{m}u^{m}\biggr)(I-Au+Qu^{2}) \\
	 &= \biggl(\sum_{m\geq 0} A_{m}u^{m}\biggr) \biggl(
	 \sum_{j=0}^{\infty}u^{2j}\biggr) (I-Au+Qu^{2}) \\
	 &= \biggl(\sum_{k\geq 0}\sum_{j=0}^{\infty} 
	 A_{k}u^{k+2j}\biggr)(I-Au+Qu^{2}) \\
	 &= \biggl(\sum_{m\geq 0}\biggl( \sum_{j=0}^{[m/2]} 
	 A_{m-2j}\biggr) u^{m}\biggr)(I-Au+Qu^{2}).
     \end{align*}     
 \end{proof}
 
 \begin{Lemma}\label{lem:eq.for.B}
     Denote by $B_{m} := A_{m} - (Q-I) \sum_{k=1}^{[m/2]}A_{m-2k}
     \in\cn(X,\G)$, for $m\geq 0$.  Then 
     
     \itm{i} $B_{0}=I$, $B_{1}=A$,
     
     \itm{ii} $B_{m} = QA_{m} - (Q-I) \sum_{k=0}^{[m/2]}A_{m-2k}$,
     
     \itm{iii} $$Tr_{\G} B_{m} = 
     \begin{cases}
	 N_{m} - Tr_{\G}(Q-I) & m \text{ even} \\
	 N_{m}  & m \text{ odd,}
     \end{cases}$$
     \itm{iv}
     $$
     \sum_{m\geq 1} B_{m}u^{m} = \left(
     Au-2Qu^{2}\right)\left(I-Au+Qu^{2}\right)^{-1}, \ \text{ for }
     |u|<\frac{1}{\a}.
     $$
 \end{Lemma}
 \begin{proof}
     $(i)$, $(ii)$ follow from computations involving bounded 
operators.
     
     $(iii)$ It follows from Lemma \ref{lem:countTail} $(ii)$ that, if
     $m$ is odd, 
     $$Tr_{\G} B_{m} = Tr_{\G}(A_{m}) - t_{m} = N_{m} ,$$ 
     whereas, if $m$ is even, 
     $$
     Tr_{\G} B_{m} = Tr_{\G}(A_{m}) - t_{m} - Tr_{\G}((Q-I)A_{0}) =
     N_{m} - Tr_{\G}(Q-I).
     $$     
     $(iv)$ 
     \begin{align*}
	 \biggl( \sum_{m\geq 0} B_{m}u^{m} \biggr)& (I-Au+Qu^{2}) 
	 \\
	 & = \biggl( Q\sum_{m\geq 0} A_{m}u^{m} - (Q-I)\sum_{m\geq
	 0}\sum_{j=0}^{[m/2]} A_{m-2j}u^{m}\biggr) (I-Au+Qu^{2}) \\
	 & = Q(1-u^{2})I - (Q-I)\biggl( \sum_{m\geq
	 0}\sum_{j=0}^{[m/2]} A_{m-2j}u^{m}\biggr) (I-Au+Qu^{2})\\
	 & = (1-u^{2})Q - (Q-I) = I-u^{2}Q,
     \end{align*}
     where the second equality follows by Lemma \ref{lem:eq.for.A}
     $(i)$ and the third equality follows by Lemma \ref{lem:eq.for.A}
     $(ii)$.  Since $B_{0}=I$, we get
     \begin{align*}
	 \biggl( \sum_{m\geq 1} B_{m}u^{m} \biggr) (I-Au+Qu^{2}) &=
	 I-u^{2}Q - B_{0}(I-Au+Qu^{2})\\
	 & = Au-2Qu^{2}.
     \end{align*}
 \end{proof}
 
 \begin{Lemma}\label{lem:Lemma3}
     Let $f:u\in B_{\eps}\equiv \{u\in\bc: |u|<\eps\} \mapsto
     f(u)\in\cn(X,\G)$, be a $C^{1}$- function, $f(0)=0$, and
     $\|f(u)\|<1$, for all $u\in B_{\eps}$.  Then 
     $$
     Tr_{\G}\left(
     -\frac{d}{du} \log(I-f(u)) \right) = Tr_{\G}\left(
     f'(u)(I-f(u))^{-1}\right).
     $$
 \end{Lemma}
 \begin{proof}
     To begin with, $-\log(I-f(u)) = \sum_{n\geq 1} \frac{1}{n}
     f(u)^{n}$, converges in operator norm, uniformly on compact
     subsets of $B_{\eps}$.  Moreover,
     $$
     \frac{d}{du} f(u)^{n} =
     \sum_{j=0}^{n-1} f(u)^{j}f'(u) f(u)^{n-j-1}.
     $$ 
     Therefore,
     $-\frac{d}{du} \log(I-f(u)) = \sum_{n\geq 1} \frac{1}{n}
     \sum_{j=0}^{n-1} f(u)^{j}f'(u) f(u)^{n-j-1}$, so that
     \begin{align*}
	 Tr_{\G}\biggl( -\frac{d}{du} \log(I-f(u)) \biggr) & =
	 \sum_{n\geq 1} \frac{1}{n} \sum_{j=0}^{n-1} Tr_{\G}\left(
	 f(u)^{j}f'(u) f(u)^{n-j-1} \right) \\
	 & = \sum_{n\geq 1} Tr_{\G}( f(u)^{n-1}f'(u) ) \\
	 & = Tr_{\G}\biggl( \sum_{n\geq 0}  f(u)^{n}f'(u) \biggr) \\
	 & = Tr_{\G}( f'(u)(I-f(u))^{-1} ),
     \end{align*}
     where we have used the fact that $Tr_{\G}$ is norm continuous.
 \end{proof}
 
 \begin{Cor} 
     $$
     Tr_{\G}\left( \sum_{m\geq 1} B_{m}u^{m} \right) = Tr_{\G}\left(
     -u\frac{d}{du} \log(I-Au+Qu^{2}) \right), \  |u|<\frac{1}{\a}.
     $$     
 \end{Cor}
 \begin{proof}
     It follows from Lemma \ref{lem:eq.for.B} $(iv)$ that
     \begin{align*}
	 Tr_{\G}\biggl( \sum_{m\geq 1} B_{m}u^{m} \biggr) &=
	 Tr_{\G}( (Au-2Qu^{2}) (I-Au+Qu^{2})^{-1} )\\	 
	 &= Tr_{\G} \Bigl( -u\frac{d}{du} \log(I-Au+Qu^{2}) \Bigr),
     \end{align*} 
     where the last equality follows from the previous lemma applied
     with $f(u) := Au-Qu^{2}$.
 \end{proof}
 
 Observe that for the $L^2$-Euler characteristic of $X$ we have
 $$
 \chi^{(2)}(X) := -\frac12 Tr_{\G}(Q-I) = |V(B)| - |E(B)| = \chi(B),
 $$
 where $\chi(B)$ is the Euler characteristic of the quotient graph
 $B=X/\G$.
 
 \begin{Thm}[Determinant formula]
     $$
     \frac{1}{Z_{X,\G}(u)} = (1-u^{2})^{-\chi(B)}
\Det_{\G}(I-Au+Qu^{2}),
     \ \text{ for } |u|<\frac{1}{\a}.
     $$ 
 \end{Thm}
 \begin{proof}
     \begin{align*}
	 Tr_{\G}\biggl( \sum_{m\geq 1} B_{m}u^{m} \biggr) &=
	 \sum_{m\geq 1} Tr_{\G}(  B_{m} ) u^{m}\\
	 &= \sum_{m\geq 1} N_{m}u^{m} - \sum_{k\geq 1} Tr_{\G}(Q-I) u^{2k} 
	 \\
	 &= \sum_{m\geq 1} N_{m}u^{m} - Tr_{\G}(Q-I) \frac{u^{2}}{1-u^{2}},
     \end{align*}
     where the second equality follows by Lemma \ref{lem:eq.for.B}
     $(iii)$.  Therefore,
     \begin{align*}
	 u\frac{d}{du} \log Z_{X,\G}(u) & = \sum_{m\geq 1} N_{m}u^{m} \\
	 &= Tr_{\G}\left( -u\frac{d}{du} \log(I-Au+Qu^{2}) \right) -
	 \frac{u}{2}\frac{d}{du} \log(1-u^{2}) Tr_{\G}(Q-I)
     \end{align*}
     so that, dividing by $u$ and integrating from $u=0$ to $u$, we
     get 
     $$
     \log Z_{X,\G}(u) = - Tr_{\G}\left( \log(I-Au+Qu^{2}) \right)
     -\frac12 Tr_{\G}(Q-I) \log(1-u^{2}),
     $$
     which implies that, for $|u|<\frac{1}{\a}$, we have
     $$
     \frac{1}{Z_{X,\G}(u)} = (1-u^{2})^{\frac12 Tr_{\G}(Q-I)}
     \cdot\exp Tr_{\G} \log(I-Au+Qu^{2}).
     $$
 \end{proof}

\section{Functional equations}

 In this final section, we obtain several functional equations for the
 Ihara zeta functions of $(q+1)-$regular graphs, $i.e.$ graphs with
 $\deg(v)=q+1$, for any $v\in VX$.  The various functional equations
 correspond to different ways of completing the zeta functions, as is
done in \cite{StTe} for finite graphs.
 
\begin{Lemma} \label{prop:holomorphy}
    Let $X$ be a $(q+1)$-regular graph and  $\D(u) := (1+qu^2)I-uA$.
Then 
    
    \itm{i} $\chi^{(2)}(X)=\chi(B)= |V(B)|(1-q)/2\in\bz$,
    
    \itm{ii} $\displaystyle Z_{X,\G}(u) = (1-u^2)^{\chi(B)}
    \Det_{\G}(\D(u))^{-1}$, for $|u| < \frac{1}{q}$,
  
    \itm{iii} by using the determinant formula in $(ii)$, $Z_{X,\G}$
    can be extended to a function holomorphic at least in the open set
    $$
    \O:=\br^2 \setminus \left(\set{(x,y)\in\br^2: x^2+y^2=\frac{1}{q}}
    \cup \set{(x,0)\in\br^2: \frac{1}{q}\leq |x|\leq 1 }\right).
    $$ 
    See figure \ref{fig:Omega}.
     \begin{figure}[ht]
	 \centering
	 \psfig{file=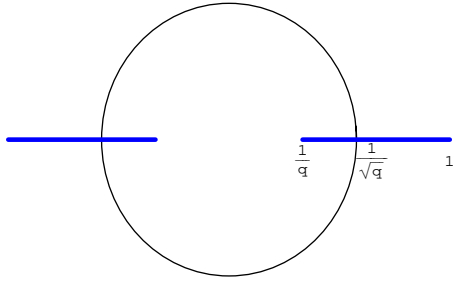,height=1.5in}
	 \caption{The open set $\Omega$}
	 \label{fig:Omega}
     \end{figure}
     
     \itm{iv} $\displaystyle \Det_\G\Bigl(\D ( \frac{1}{qu})\Bigr) =
     (qu^2)^{-|VB|} \Det_\G(\D(u))$, for $u\in\O\setminus \set{0}$.
     
\end{Lemma}
\begin{proof}
    $(i)$ This follows by a simple computation.
    
    $(ii)$ This follows from $(i)$.

    $(iii)$ Let us observe that 
    $$
    \s(\D(u)) = \set{1+qu^2-u\l: \l\in\s(A)} \subset \set{1+qu^2-u\l:
    \l\in[-d,d]}.
    $$
    It follows that $0\not\in\conv\s(\D(u))$ at least for $u\in\bc$
    such that $1+qu^2-u\l\neq0$ for $\l\in[-d,d]$, that is for $u=0$
    or $\frac{1+qu^2}{u}\not\in[-d,d]$, or equivalently, at least for
    $u\in\O$.  The rest of the proof follows from Corollary
    \ref{cor:det.analytic}.
    
    $(iv)$ This follows by Proposition \ref{properties} $(i)$ and the
    fact that $Tr_\G(I_V) = |VB|$.
\end{proof}

  \begin{Prop} [Functional equations]
     Let $X$ be $(q+1)$-regular. Then, for all $u\in\O$, we have

     \itm{i} $\La_{X,\G}(u) := (1-u^{2})^{-\chi(B)}(1-u^{2})^{|VB|/2}
     (1-q^{2}u^{2})^{|VB|/2} Z_{X,\G}(u) =
     -\La_{X,\G}\Bigl(\frac{1}{qu}\Bigr)$,
     
     \itm{ii} $\xi_{X,\G}(u) := (1-u^{2})^{-\chi(B)} (1-u)^{|VB|}
     (1-qu)^{|VB|} Z_{X,\G}(u) = \xi_{X,\G}\Bigl(\frac{1}{qu}\Bigr)$,
     
     \itm{iii} $\Xi_{X,\G}(u) := (1-u^{2})^{-\chi(B)}
     (1+qu^{2})^{|VB|} Z_{X,\G}(u) =
     \Xi_{X,\G}\Bigl(\frac{1}{qu}\Bigr)$.
 \end{Prop}
 \begin{proof}
     $(i)$ 
     \begin{align*}
	 \La_{X}(u) & = (1-u^{2})^{|VB|/2} (1-q^{2}u^{2})^{|VB|/2}
	 \Det_{\G}(\D(u))^{-1} \\
	 &= u^{|VB|}\Bigl(\frac{q^{2}}{q^{2}u^{2}}-1 \Bigr)^{|VB|/2}
	 (qu)^{|VB|} \Bigl(\frac{1}{q^{2}u^{2}}-1 \Bigr)^{|VB|/2}
	 \frac{1}{(qu^{2})^{|VB|}}\Det_{\G}\Bigl( \D(\frac{1}{qu})
	 \Bigr)^{-1}\\
	 &= -\La_{X}\Bigl(\frac{1}{qu}\Bigr).
     \end{align*}
     $(ii)$
     \begin{align*}
	 \xi_{X}(u) & = (1-u)^{|VB|} (1-qu)^{|VB|} \Det_{\G}(\D(u))^{-1}\\
	 & = u^{|VB|} \Bigl( \frac{q}{qu} -1 \Bigr)^{|VB|} (qu)^{|VB|}
	 \Bigl( \frac{1}{qu} -1 \Bigr)^{|VB|}
	 \frac{1}{(qu^{2})^{|VB|}}\Det_{\G}\Bigl( \D(\frac{1}{qu})
	 \Bigr)^{-1}\\
	 & = \xi_{X}\Bigl(\frac{1}{qu}\Bigr).
     \end{align*}
     $(iii)$
     \begin{align*}
	 \Xi_{X}(u) & = (1+qu^{2})^{|VB|} \Det_{\G}(\D(u))^{-1} \\
	 &= (qu^{2})^{|VB|} \Bigl( \frac{q}{q^{2}u^{2}} +1 \Bigr)^{|VB|}
	 \frac{1}{(qu^{2})^{|VB|}}\Det_{\G}\Bigl(
	 \D(\frac{1}{qu}) \Bigr)^{-1}\\
	 &= \Xi_{X}\Bigl(\frac{1}{qu}\Bigr).
     \end{align*}
 \end{proof} 
 
  \begin{ack}
     The second and third named authors would like to thank
     respectively the University of California, Riverside, and the
     University of Roma ``Tor Vergata'' for their hospitality at
     different stages of the preparation of this paper.
 \end{ack}


\end{document}